\documentclass[12pt]{amsart}

\usepackage[margin=1in]{geometry}

\usepackage{amsmath,amscd,amssymb,amsfonts,latexsym,wasysym, mathrsfs, mathtools,hhline,color}
\usepackage[all, cmtip]{xy}
\usepackage{csquotes}
\usepackage{url}
\usepackage{comment}

\definecolor{hot}{RGB}{65,105,225}

\usepackage[pagebackref=true,colorlinks=true, linkcolor=hot ,  citecolor=hot, urlcolor=hot]{hyperref}
\usepackage{ textcomp }
\usepackage{ tipa }
\usepackage{graphicx,enumerate}

\usepackage{geometry}
\theoremstyle{plain}
\newtheorem{theorem}{Theorem}[section]
\newtheorem{prop}[theorem]{Proposition}

\newtheorem{lm}[theorem]{Lemma}

\newtheorem{cor}[theorem]{Corollary}

\newtheorem{thrm}[theorem]{Theorem}

\theoremstyle{definition}

\newtheorem{defn}[theorem]{Definition}

\newtheorem{rmk}[theorem]{Remark}

\newtheorem{ex}[theorem]{Example}
\newtheorem*{ex*}{Example}

\def\be{\begin{equation}}
\def\ee{\end{equation}}

\def\bt{\begin{thrm}}
\def\et{\end{thrm}}

\def\bc{\begin{cor}}
\def\ec{\end{cor}}

\def\br{\begin{rmk}}
\def\er{\end{rmk}}

\def\bp{\begin{prop}}
\def\ep{\end{prop}}

\def\bl{\begin{lm}}
\def\el{\end{lm}}

\def\bex{\begin{ex}}
\def\eex{\end{ex}}

\def\bd{\begin{defn}}
\def\ed{\end{defn}}

\newcommand{\C}{\mathbb{C}}
\newcommand{\Z}{\mathbb{Z}}
\newcommand{\Hom}{\mathrm{Hom}}

\newcommand{\K}{\mathbb{K}}

\newcommand{\cha}{\mathrm{char}}

\newcommand{\orb}{\mathrm{orb}}
\newcommand{\M}{\mathbf{M}}

\newcommand{\im}{\mathrm{im}}
\newcommand{\homo}{\mathrm{Hom}}
\newcommand{\sV}{\mathcal{V}}

\newcommand{\sA}{\mathcal{A}}

\newcommand{\rank}{\mathrm{rank }}

\newcommand{\bm}[1]{\mbox{\boldmath{$#1$}}}

\title[]{The homology growth for finite abelian covers of smooth quasi-projective varieties}

\author{Fenglin Li}
\address{Fenglin Li: School of Mathematical Sciences, University of Science and Technology of China,  96 Jinzhai Road, Hefei Anhui 230026 China}
\email{fl0125@mail.ustc.edu.cn}

\author{Yongqiang Liu}
\address{Yongqiang Liu: The Institute of Geometry and Physics, University of Science and Technology of China, 96 Jinzhai Road, Hefei Anhui 230026 China} 
\email{liuyq@ustc.edu.cn}

\date{\today}



\begin{document}
\maketitle
\begin{abstract}
Let $X$ be a complex smooth quasi-projective variety with a fixed epimorphism $\nu\colon\pi_1(X)\twoheadrightarrow H$, where $H$ is a finitely generated abelian group with $\mathrm{rank}H\geq 1$. In this paper, we study the asymptotic behaviour of Betti numbers with all possible field coefficients and the order of the torsion subgroup of singular homology associated to $\nu$, known as the $L^2$-type invariants. When $\nu$ is orbifold effective, we give explicit formulas of these invariants at degree 1. This generalizes the authors' previous work for $H\cong \Z$.
\end{abstract}
\section{Introduction}

There is a general principle to consider a classical invariant of a finite CW complex $X$ and to define its analogues for some  covering space of $X$. This leads to the $L^2$-type invariants.  Atiyah  \cite{Ati76} introduced
the notion of $L^2$-Betti numbers in the context of a regular covering  of a closed
Riemannian manifold. After that,  there have been vast literatures for the $L^2$-invariant theory, see \cite{Luc02}.   A particular important result is L\"uck's approximation theorem \cite{Luc94}, which states  that the $L^2$-Betti numbers of the
universal cover of a finite CW complex can be found as limits of normalized Betti numbers
of finitely sheeted normal coverings.  Then it becomes a classical subject to study $L^2$-type invariants focusing on its approximation by towers of finite coverings, see e.g. \cite{Luc13,Luc16}.

\bigskip

Let $X$ be a connected finite CW complex with a fixed epimorphism $\nu\colon\pi_1(X)\twoheadrightarrow H$, where $H$ is a finitely generated abelian group with $\rank H=n \geq 1$.
We fix an isomorphism $H\cong \Z^n\oplus T$, where $T$ is a finite abelian group.  Consider $\Z^n$ as a subgroup of $H$ under this isomorphism.  For a subgroup $\Gamma \subset \Z^n$ of finite index, we set
$$\langle\Gamma\rangle=\min \left\{ \sum_{j=1}^n x_j^2  \mid x=(x_1,\ldots,x_n)\in \Gamma, x\neq 0 \right\} $$
Let $X^{\nu,\Gamma}$ denote the covering space of $X$ associated to the corresponding composition of $ \nu$ and the quotient map $H\to H/\Gamma$.
 Consider the following limits:
\[ \alpha_i(X^{\nu},\K)\coloneqq\lim_{\langle\Gamma\rangle\rightarrow\infty}\frac{\dim H_{i}(X^{\nu,\Gamma},\mathbb{K})}{|H/\Gamma|} \]
for any field coefficients $\K$ and
\[ \M_i(X^{\nu})\coloneqq\limsup_{\langle\Gamma\rangle\rightarrow\infty}\frac{\log|H_{i}(X^{\nu,\Gamma},\mathbb{Z})_{\mathrm{tor}}|}{|H/\Gamma|}. \]
Here $ |H_{i}(X^{\nu,\Gamma},\mathbb{Z})_{\mathrm{tor}}|$ denotes the order of the torsion part of $H_{i}(X^{\nu,\Gamma},\mathbb{Z}).$ These two limits are particular cases of $L^2$-type invariants. Such kind of $L^2$-type invariants have been studied by many people (see \cite{Luc13,Luc16} for more results in this direction).

\medskip

When $H$ is free abelian,  Linnell, L\"{u}ck, Sauer \cite[Theorem 0.2]{Luc11} (for $\cha(\K)=0$) and Ab\'ert, Jaikin-Zapirain, Nikolov \cite[Theorem 17]{AJN11} (for $\cha(\K)>0$) showed  that the first limit always exists. Meanwhile, the second limit also exists and can be computed by the Mahler measure of the Alexander polynomial, see \cite{Le14}. 

Moreover, if $H$ has non-trivial  torsion part (i.e., $T\neq 0$), the computation of these limits can be reduced  to the free abelian case using a finite cover trick, see section 3.

When  $X$ is a complex smooth quasi-projective variety and $H\cong \Z$, the authors \cite{LL21} gave concrete formulas for these limits (at degree 1) in terms of the geometric information of $X$. The main results in this paper are to generalize these formulas to arbitrary finitely generated abelian group $H$ with $\rank H\geq 1$,
when $\nu$ is an orbifold effective morphism.

\medskip

We first give the definition of orbifold morphism.

\bd Let $X$ be a smooth complex quasi-projective variety. An algebraic map $f\colon X \to \Sigma_{g,r}$  is called an orbifold map, if $f$ is surjective, has connected generic fiber and $\Sigma_{g,r}$ is a smooth algebraic curve of genus $g$ with $r$ points removed. We always assume that $\Sigma_{g,r} \neq \mathbb{CP}^1, \C^1$. There exists a maximal Zariski open subset $U\subset \Sigma_{g,r}$ such that $f$ is a fibration over $U$. Say $B=\Sigma_{g,r}-U$ (could be empty) has $s $ points, denoted by $\{q_1,\ldots,q_s\}$.  We assign
the multiplicity $\mu_j$ (the $\gcd$ of the coefficients of the divisor $f^* q_j$) of the fiber   $f^{*}(q_j)$  to the point $q_j$.   
Such orbifold map $f$ is called of type $(g,r,{\bm \mu})$, 
where ${\bm \mu}=(\mu_1,\ldots,\mu_s)$. When $B=\emptyset$, $\prod_{j=1}^s \mu_j=1$ by convention. \ed

The orbiford group $\pi_1^{\orb}(\Sigma_{g,r},{\bm \mu})$ associated to these data is  defined as
\[ \pi_{1}^{\orb}(\Sigma_{g,r},{\bm \mu})\coloneqq \pi_{1}(\Sigma_{g,r}\backslash \{q_{1}, \ldots, q_{s}\})/ \langle \gamma_{j}^{\mu_{j}}=1 \text{ for all } 1\leq j \leq s\rangle, \]
where $\gamma_{j}$ is a meridian of $q_{j}$.
An orbifold map  $f\colon X \to \Sigma$ of type $(g,r,{\bm \mu})$ induces an surjective map to the orbifold group (see \cite[Proposition 1.4]{ACM10})
$$ f_* \colon \pi_1(X) \twoheadrightarrow \pi_1^{\orb}(\Sigma_{g,r},{\bm \mu}).$$
When $\Sigma_{g,r}$ is clear in the context, we simply write $\Sigma$.

\bd \label{def orbifold effective}
Let $X$ be a smooth complex quasi-projective variety with an epimorphism $\nu\colon \pi_1(X)\twoheadrightarrow H$. We say that $\nu$ is orbifold effective if there is an orbifold map $f\colon X\rightarrow\Sigma_{g,r}$  such that $\nu$ factors through $f_*$ as follows:
\[\xymatrix{
\pi_1(X)\ar@{->>}"1,3"^{\nu}\ar@{->>}[dr]_{f_*} & & H \\
 & \pi_1^{\orb}(\Sigma_{g,r}, {\bm \mu})\ar@{->>}[ur] &
}.\]
 We  say that $\nu$ is  orbifold effective by $f$ and call $\nu$ being of type $(g,r,{\bm \mu})$.
\ed

Our main result is the following:

\bt\label{main thm}
Let $X$ be a complex smooth quasi-projective variety with a fixed epimorphism $\nu\colon\pi_1(X)\twoheadrightarrow H$, where $H$ is a finitely generated abelian group with $\mathrm{rank}H \geq 1$. Suppose that $\nu$ is  orbifold effective of type $(g,r, {\bm \mu})$. Let $\K$ be a field with $\cha(\K)=p\geq 0$. 
\begin{itemize}
\item[(a)] If $H$ is a free abelian group, then
$$  \alpha_1(X^\nu,\K)=  2g+r-2+ \#\{ j  \mid p\ \mathrm{ divides }\ \mu_j\} $$
and
$$  \M_1(X^\nu)=\sum_{j=1}^s \log \mu_j. $$
\item[(b)] If $H$ has non-trivial torsion part, then we have
$$
\alpha_1(X^{\nu},\K)=2g+r-2+\sum_{j=1}^s\left(1-\frac{1}{m_j}\right)+\sum_{1\leq j\leq s,  p\mid\frac{\mu_j}{m_j}}\frac{1}{m_j}.
$$
and
$$
\M_1(X^{\nu},\K)=\sum_{j=1}^s\frac{1}{m_j}\log\frac{\mu_j}{m_j},
$$
where $m_j$ is a positive integer dividing $\mu_j$ and it only depends on $X$ and $\nu$. For its definition, see Remark \ref{rem how to define m}.
\end{itemize}
\et

The proof of Theorem \ref{main thm} is based on the theory of cohomogy jump loci.
In section 2, we recall some basic properties of cohomology jump loci and give the proof of Theorem  \ref{main thm} for the free abelian case.  Section 3 is devoted to the proof of Theorem  \ref{main thm} for the non-free abelian case using a finite cover trick.

\bigskip

\textbf{Acknowledgments.} 
The second named author  is partially supported by National Key Research and Development Project SQ2020YFA070080, NSFC-12001511, the Project of Stable Support for Youth Team in Basic Research Field CAS (YSBR-001), a research fund from University of Science and Technology of China,  the project “Analysis and Geometry on Bundles” of the Ministry of Science and Technology of the People's Republic of China and  Fundamental Research Funds for the Central Universities.

\section{The Free Abelian Case}
In this section we consider the case when $H$ is a free abelian group with rank $n$, say $\Z^n$.
\subsection{Limits of Betti numbers}

We recall the definition of cohomology jump loci. Let $X$ be a connected finite CW complex with $\pi_1(X)=G$. The group of $\K$-valued  characters, $ \homo(G,\K^*)$, is a commutative affine algebraic group. Each character $\rho \in \homo(G,\K^*)$ defines a rank one local system on $X$, denoted by $L_{\rho}$. Note that $\mathrm{Hom}(G,\K^*)$ only depends on $H_1(X, \Z)$, the abelianization of $G$.
\bd
The cohomology jump loci of $X$ are defined as
$$\sV^i_k(X,\K)\coloneqq \lbrace \rho\in \homo(G,\K^*) \mid \dim_{\K} H^{i}(X, L_{\rho})\geq k \rbrace.$$
When $k=1$, we simply write $\sV^i(X,\K)$.
\ed

Cohomology jump loci  are closed sub-varieties of $\homo(G,\K^*)$ and homotopy invariants of $X$.
 In degree 1, $\sV^1_k(X,\K)$  depends only on $\pi_1(X)$ (e.g. see \cite[Section 2.2]{Suc11}).

The map $\nu\colon G\twoheadrightarrow H\cong \Z^n$ induces an embedding $(\K^*)^n\subset\homo(G,\K^*)$. For a tuple $\lambda=(\lambda_1,\ldots,\lambda_n)\in(\K^*)^n$, let $\nu^{-1} L_{\lambda}$ denote the corresponding rank one local system on $X$ whose monodromy representation factors through $\nu$.

\bp  \label{prop betti number} Let $\K$ be an algebraically closed field.
With the  notations  as above,  for any $i\geq 0$ and $\lambda \in (\K^*)^n$ being general we have
$$ \alpha_i(X^\nu,\K)=\dim H^i(X, \nu^{-1} L_\lambda).$$
In particular, $\alpha_i(X^\nu,\K)$ is always an integer.
\ep
\begin{proof}
Let $\pi \colon X^{\nu,\Gamma}\rightarrow X$ denote the covering map. Then
$$H^{i}(X^{\nu,\Gamma},\K)=H^{i}(X, \pi_*\K),$$
where $\pi_*\K$ is the  push forward of the $\K$-constant sheaf on $X^{\nu,\Gamma}$, hence a rank $|\Z^n/\Gamma|$ local system.

\cite[Theorem 0.2]{Luc11} and \cite[Theorem 17]{AJN11} show that the limit $ \alpha_i(X^\nu,\K)$ always exists.  By choosing a sub-sequence,  we may assume that $\cha(\K)\nmid \vert \Z^n/\Gamma \vert$. Since $\Gamma$ is a sub-group of $\Z^n$ with finite index, $\Gamma$ is also a free abelian group with rank $n$.  Using Smith normal form, without loss of generality we say $\Gamma= N_1 \Z \oplus \cdots \oplus N_n \Z $ for  a tuple of positive integers where $N_n \mid N_{n-1} \mid \cdots \mid N_1$. In particular, $\cha(\K)\nmid N_1$. Then the local system $\pi_*\K$ decomposes as the direct sum of $\vert \Z^n/\Gamma \vert$-many rank one local systems.   
Say  $L_{\lambda}$ with $\lambda=(\lambda_1,\ldots,\lambda_n)\in(\K^*)^n$ is one of the direct sum factors.
Then  $\lambda_j^{N_{j}}=1$ for any $1\leq j\leq n$.  Suppose that $a$ is the number such that
\[ (\K^*)^n\subseteq\sV^i_a(X,\K)\ \mathrm{and}\ (\K^*)^n\nsubseteq\sV^i_{a+1}(X,\K). \]
Note that $(\K^*)^n\cap\sV^i_{a+1}(X,\K)$ is a subvariety of $(\K^{*})^n$ with less dimension. There exists a hypersurface $V$ such that
\[ (\K^*)^n\cap\sV^i_{a+1}(X,\K)\subset V=\{ (t_1,\ldots,t_n)\in(\K^*)^n | u(t_1,\ldots,t_n)=0 \}, \]
where $u$ is a polynomial. Without loss of generality, suppose the degree of $t_1$ for $u$ is $d\geq 1$. Then for a fixed $(\lambda_2,\ldots,\lambda_n)$, $u(t_1,\lambda_2,\ldots,\lambda_n)$ has at most $d$-many solutions. Hence we have
\[ a \cdot \prod_{j=1}^n N_j\leq \dim H^{i}(X^{\nu,\Gamma},\mathbb{K})\leq a\cdot \prod_{j=1}^n N_j+ c\cdot d \cdot \prod_{j=2}^n N_j ,\]
where $c$ is some constant number which only depends on $X$ (e.g. $c$ can be taken as the number of $i$-cells in $X$).
It implies \[ a\leq \frac{\dim H^{i}(X^{\nu,\Gamma},\mathbb{K})}{\prod_{j=1}^n N_j}\leq a+c\cdot \frac{d}{N_1}  ,\]
For field coefficients, $ \dim H^{i}(X^{\nu,\Gamma},\mathbb{K})=\dim H_{i}(X^{\nu,\Gamma},\mathbb{K}) $.
Taking $\langle \Gamma \rangle\rightarrow\infty$, we are done since $N_1$ goes to infinity.
\end{proof}

When $\K=\C$, the cohomology jump loci of complex smooth quasi-projective variety  have been  intensively studied. In particular, the following structure theorem for $\sV^i_k(X,\C)$ put strong constraints for the homotopy type of complex smooth quasi-projective variety. It is contributed by many people and we name a few here: Green-Lazarsfeld \cite{GL91}, Simpson \cite{Sim93}, Arapura \cite{Ara97}, Dimca-Papadima \cite{DP14}, Dimca-Papadima-Suciu \cite{DPS}, etc. It is finalized by Budur and Wang in \cite{BW15,BW20}.
\bt \label{structure theorem} \cite{BW15,BW20} If $X$ is a complex smooth variety, then $ \sV^i_k(X,\C)$ is a finite union of torsion translated sub-tori of $\Hom(G,\C^*)$.
\et

\subsection{Limits of torsion}

The second type limit also exists, see \cite[Theorem 5]{Le14}. It can be computed by the Mahler measure of the $i$-th integral Alexander polynomial $\Delta_i(X^{\nu})\in\Z[t^{\pm}_1,\ldots,t^{\pm}_n]$. Let us recall the definitions of multivariable Alexander polynomials and Mahler measures.

Recall that $X$ is a connected finite CW-complex with a group epimorphism $\nu\colon\pi_{1}(X)\twoheadrightarrow\mathbb{Z}^n$.
Then the group of covering transformations of the covering space  $X^{\nu}$ is isomorphic to $\Z^n$ and acts on it.   By choosing fixed lifts of the cells of $X$ to $X^{\nu}$, we obtain a free basis for the cellular chain complex  of $ X^\nu$ as   $R_n$-modules, where $R_n= \Z[\Z^n]=\Z[t^{\pm}_1,\ldots,t^{\pm}_n]$.   So the cellular chain complex of $X^{\nu}$, $C_{*}(X^{\nu}, \Z)$, is a bounded complex of finitely generated free $R_n$-modules:
\be \label{chain compelx}
 \cdots  \to  C_{i+1}(X^\nu, \Z) \overset{\partial_{i}}{\to} C_i(X^\nu, \Z) \overset{\partial_{i-1}}{\to} C_{i-1}(X^\nu, \Z)  \overset{\partial_{i-2}}{\to}   \cdots \overset{\partial_0}{\to} C_0(X^\nu, \Z)  \to 0 .
\ee
With the above free basis for $C_*(X^\nu,\Z)$, $\partial_i$ can be written down as a matrix with entries in $R_n$. Note that $R_n$ is a Notherian UFD.
Assume that $\partial_i$ has rank $r_i$.
Let $\Delta_i(X^\nu)$ denote the greatest common divisor of all non-zero $(r_i\times r_i)$-minors of $\partial_i$.
When $\partial_i=0$, $\Delta_i(X^\nu)=1$ by convention. 

\bd \label{def Alexander}$ \Delta_i(X^\nu)$  is called {\it the $i$-th $n$-variable Alexander polynomial of $(X,\nu)$}. Then $ \Delta_i(X^\nu)$ is defined uniquely up to a multiplication with a unit of $R_n$.
\ed

Let $h\in R_n$ be a nonzero polynomial. The Mahler measure of $h$ is defined by
\[ \M(h)\coloneqq\int_{(S^1)^n} \log |h(s)|\mathrm{d}s,\]
where $\mathrm{d}s$ indicates integration with respect to normalized Haar measure, and $(S^1)^n$ is the multiplicative
subgroup of $n$-dimensional complex space $\C^n$ consisting of all vectors $(s_1, \ldots,s_n)$ with $|s_1|= \cdots =|s_n|=1$. Here $h$ is regarded as a function on $\C^n$.\par

Recall that an element $h\in R_n$ is called a generalized cyclotomic polynomial \cite[page 47]{Sch} if it is of the form $h(t_1,\ldots,t_n)=t^{\mathbf{m}}\Phi(t^{\mathbf{n}})$, where $\mathbf{m},\mathbf{n}\in\Z^n,\mathbf{n}\neq 0$, and $\Phi$ is a cyclotomic polynomial in a single variable.

\bt\cite[Theorem 19.5]{Sch}\label{general cyclotomic}
Let $h\in R_n$. Then $\M(h)=0$ if and only if $\pm h$ is a product of generalized cyclotomic polynomials.
\et

\bt\cite[Theorem 5]{Le14}\label{Le thm}
With the notations above, we have
\[ \limsup_{\langle \Gamma\rangle\rightarrow\infty}\frac{\log|H_{i}(X^{\nu,\Gamma},\mathbb{Z})_{\mathrm{tor}}|}{|\Z^n/\Gamma|}=\M(\Delta_i(X^{\nu})). \]
When $n=1$, then $\limsup$ can be replaced by the ordinary $\lim$.
\et

  Theorem \ref{structure theorem} implies the following property for the Alexander polynomial associated to the pair $(X,\nu)$.  For the one-variable version of the following result, see \cite[Proposition 1.4]{BLW} and \cite[Proposition 3.7]{LL21}.

\bp \label{lim tor}
Let $X$ be a complex smooth variety with a fixed epimorphism $\nu\colon\pi_1(X)\twoheadrightarrow \Z^n$. Then $\Delta_i(X^{\nu})$ is a product of an positive integer $c_i$ with some generalized cyclotomic polynomials, where $c_i$ is the leading coefficient of $\Delta_i(X^{\nu}) $. In particular,  $\M(\Delta_i(X^{\nu}))=\log c_i$.
\ep
\begin{proof}
   The structure theorem for cohomology jump loci $\sV^i_k(X,\C)$ implies that $(\C^*)^n \cap \sV^i_k(X,\C)$ is a finite union of torsion translated sub-tori.
  From now on, we consider  $\Delta_i(X^\nu)$ as an element in $\C[t^\pm_1,\ldots,t^{\pm}_n]$.
We claim that the irreducible factors of  $\Delta_i(X^\nu)$ can be observed by the irreducible hyersurfaces in $(\C^*)^n \cap \sV^i_k(X,\C)$.
 Say $h$ is a irreducible factor of $\Delta_i(X^\nu)$. Then $h$ generates a prime ideal with height 1. Let $\C[t^\pm_1,\ldots,t^{\pm}_n]_{(h)}$ be its localization at the prime ideal  $(h)$, which is a PID. Consider $ H_i(X^\nu ,\C)$ as a finitely generated $\C[t^\pm_1,\ldots,t^{\pm}_n]$-module and its localization $ H_i(X^\nu ,\C)_{(h)}$.
 Then by a similar proof to \cite[Theorem 4.2]{DN},
 we have
 $$ \dim H_i(X, \nu^{-1} L_\lambda)= \rank  H_i(X^\nu ,\C)_{(h)} + J_i +J_{i-1}, $$
 where $\lambda $ is a general point in the zero locus of $h$, $ \rank  H_i(X^\nu ,\C)_{(h)}$ is its rank as a finitely generated module over the PID $\C[t^\pm_1,\ldots,t^{\pm}_n]_{(h)}$ and $J_i$ is the number of direct summands of the torsion part of  $H_i(X^\nu ,\C)_{(h)}$.
 By Universal Coefficients Theorem over the PID $\C[t^\pm_1,\ldots,t^{\pm}_n]_{(h)}$,
one can translate from homology to cohomology. Then by structure theorem, we get that  $h$ has the form $\prod_{j=1}^n t_j^{d_j}= \epsilon$, where $(d_1,\ldots,d_n)$ is a non-zero tuple of integers and $\epsilon$ is a roots of unity.  Since $\Delta_i(X^\nu)$ has integer coefficients,  $\Delta_i(X^\nu)$ is a product of generalized cyclotomic polynomial with some integer $c_i$ by Galois theory.
Then the claim follows. In particular, it follows from Theorem \ref{general cyclotomic} that $\M(\Delta_i(X^{\nu}))=\log c_i$.
 \end{proof}

\begin{proof}[Proof of Theorem \ref{main thm}(a)]

Note that since $b_i(X,\K)$ only depend on  $\cha(\K)$, not on the specific choice of the field $\K$. So without loss of generality, we can assume that $\K$ is algebraically closed. The first equation  follows from Proposition \ref{prop betti number} and \cite[Theorem 1.7]{LL21}.

For the second equation, consider the following commutative diagram
\[\xymatrix{
\pi_1(X)\ar@{->>}"1,3"^{\nu}\ar@{->>}[dr]_{f_*} & & \Z^n\ar@{->>}^{\kappa}[r] & \Z, \\
 & \pi_1^\orb(\Sigma_{g,r})\ar@{->>}[ur] &
}\]
where $\kappa$ is a surjective morphism to $\Z$.
Say $\kappa$ is represented by a tuple $(a_1,\ldots,a_n)$ with $\gcd(a_1,\ldots,a_n)=1$. Then it is easy to see if we substitute $(t_1,\ldots,t_n)$ by $(t^{a_1},\ldots,t^{a_n})$ in the matrix $\partial_1$ of the chain complex (\ref{chain compelx}), we get the corresponding matrix for the chain complex of $X^{\kappa\circ \nu}$. By choosing a general tuple $(a_1,\ldots,a_n)$, the rank of $\partial_1$ does not change after substituting. Hence $\Delta_1(X^{\kappa\circ\nu})=\Delta_1(X^{\nu})(t^{a_1},\ldots,t^{a_n})$ in this case. Note that $\kappa\circ \nu$ is an orbifold effective morphism of the type $(g,r,{\bm \mu})$. By \cite[Theorem 1.11]{LL21} we know that $\Delta_1(X^{\kappa\circ\nu})$ has leading coefficient $\prod_{j=1}^s\mu_j$. It coincides with the leading coefficient of $\Delta_1(X^{\nu})$. Then the claim follows from the above proposition.
\end{proof}

\section{The General Case}
In this section, we study the case where $H$ has non-trivial torsion part.
Say $H\cong \Z^n \oplus T$, where $T$ is a finitely generated torsion abelian group. Let $q$ be the quotient map $H\rightarrow T$ and $X^T$ denote the corresponding finite cover  of $X$ associated to $q\circ\nu$.
Then it is easy to see that there exists an epimorphism $\nu^T \colon \pi_1 (X^T) \to \Z^n$ such that the following  diagram commutes
\[ \xymatrix{
\pi_1(X^T)\ar@{>>}^{\nu^T}[r]\ar[d] & \Z^n\ar@{>>}[r]\ar[d] & \Z^n/\Gamma\ar[d]\\
\pi_1(X)\ar@{>>}^{\nu}[r] & \Z^n\oplus T\ar@{>>}[r] & \Z^n/\Gamma\oplus T.
}\]
Due to the choice of $\Gamma$, we can identify the finite index covering spaces associated to the two composed horizontal maps. Then we have the following result, which reduces the computations to the free abelian case.
\bl \label{lem finite cover}
With the notations above, we have the following equations
\be
\alpha_1(X^{\nu},\K)=\frac{1}{|T|}\alpha_1((X^T)^{\nu^T},\K)
\ee
and
\be
\M_1(X^{\nu})=\frac{1}{|T|}\M_1((X^T)^{\nu^T}).
\ee
\el

We can choose $X^T$ to be a smooth quasi-projective variety and the corresponding cover map $X^T \to X$ to be algebraic. 
 Assume that the epimorphism $\nu\colon\pi_1(X)\twoheadrightarrow H$ is orbifold effective by some orbifold map $f\colon X\rightarrow C$ for some smooth curve $C$.  By  projectivising, we get a map $h \colon \overline{X^T} \to \overline{C}$ for the map $X^T \overset{f\circ \pi}{\longrightarrow} C$. Using Stein factorization,  we get the following commutative diagram $$\xymatrix{
X^T \ar[d]^{g}   \ar[r] & \overline{X^T} \ar[d]^{ h'} \ar[rd]^{h} &  \\
S \ar[r]                & \overline{S}  \ar[r]^{h''} & \overline{C},
}$$
where $h''$ is a finite map, $h'$ has connected fiber and  $g\coloneqq h'\vert_{X^T}\colon  X^T  \to S\coloneqq \im(g)$. Then \cite[Lemma 2.2]{Dim07} shows that $g$ has connected generic fiber. Hence $g$ is an orbifold map and we have the following commutative diagram
\[ \xymatrix{
  X^{T}\ar^{g}[r]\ar_{\pi}[d] & S\ar^{\pi^{\prime}}[d] \\
 X\ar^{f}[r] & C,
}\]
where $\pi$ is the covering map and $\pi'$ is obtained from $h''$ by taking restrictions over $S$.

The image of the following composed map
$$ \pi_1(X^T)\to \pi_1(X)\to \pi_1^\orb(C)\to \Z^n\oplus T$$
is $\Z^n$ (the last map exists since $f$ is orbifold effective), hence 
one has the following commutative diagram \[ \xymatrix{
 \pi_1( X^{T})\ar[d]\ar[r] & \pi_1(S)\ar[d]  &  \\
 \pi_1(X)\ar[r] & \pi_1(C) \ar[r]& \Z^n.
}\] In particular, the orbifold map $g$ makes $\pi_1(X^T)\overset{ \nu^T}{\twoheadrightarrow} \Z^n$ orbifold effective.

Using  Lemma \ref{lem finite cover}, one can reduce the computation to the free abelian case. All we need to know is the orbifold information of $g$. For any point $b\in C$, let $D_b$ be a  small enough open disc of $b$. Set $T_b=f^{-1}(D_b)$.  
Next we study the point $b$ in two cases.

\begin{enumerate}
\item[Case 1:]{$b\notin B$.} Then  $F=f^{-1}(b)$ is the generic fiber of $f$ and $T_b$ is smooth over $D_b$, hence admits a trivial fibaration over $D_b$. Dimca proved the following short  exact sequence  (see \cite[Section 5]{Dim07} and \cite[Section 6.3]{Dim17})
$$H_1(F,\Z) \to H_1(X,\Z) \to H_1(\pi^\orb_1(C)) ,$$
 where $H_1(\pi_1(C)^\orb) $ is the abelization of $\pi^\orb_1(C)$.  Then the following composed group homomorphism is trivial
$$ \pi_1(F) \to \pi_1(X) \twoheadrightarrow \pi^\orb_1(C) \twoheadrightarrow H \twoheadrightarrow T,$$
since $T$ is abelian. Hence $\pi^{-1}(F)$ has $\vert T \vert$-many disjoint copies of $F$.
Note that $g$ is an orbifold map, which has connected generic fiber. Therefore $(\pi')^{-1}(b)= g (\pi^{-1}(F)) $ consists of $\vert T \vert$-many points.  In particular, $F$ is also a generic fiber for $g$.

\item[Case 2:]{$b\in B$.} 
Consider the following composed group homomorphism
$$ \pi_1(T_b)\to \pi_1(X) \twoheadrightarrow \pi_1^\orb(C)  \twoheadrightarrow H \twoheadrightarrow T.$$
Let $m_b $ denote the order of the image group of the above composed group homomorphism.
Then $\pi^{-1}(T_b)$ has $\frac{\vert T\vert}{m_b}$-many connected components. Note that $g$ is an orbifold morphism, which has connected generic fiber.  Hence  $g(\pi^{-1}(T_b))=(\pi')^{-1}(D_b)$ also has $\frac{\vert T\vert}{m_b}$-many connected components, i.e., ${\pi'}^{-1}(b)$ has $\frac{\vert T\vert}{m_b}$ many points. 
Therefore we have the following commutative diagram
\[ \xymatrix{
\pi^{-1}(T_b)\ar^{\mu_b/m_b}_g[r]\ar[d] & (\pi')^{-1}(D_b)\ar[d]_{\pi'}^{m_b} \\
T_b\ar^{\mu_b}_f[r] & D_b.
}\]
Let us explain the integers appearing in the diagram. $\mu_b$ for the bottom horizontal map means that $f$ has multiplicity $\mu_b$ over $b$.
Note that for any point $a\in D_b \setminus \{b\}$, $(\pi')^{-1}(a)$ has $\vert T \vert$-many points and $(\pi')^{-1}(b)$ has $\frac{\vert T\vert}{m_b}$-many points. Due to the construction of $g$, $(\pi')^*b$ as a divisor has the same coefficients for every point in $(\pi')^{-1}(b)$ (since composing with the deck transformation over $X^T$ does not change $g$). So $\pi'$ has multiplicity $m_b$ over  $b$. 
On the other hand, since $\pi$ is a finite cover map, the left vertical map has multiplicity 1.
Putting these together, we get that $g$ has multiplicity $ \frac{\mu_b}{m_b}$ over every point in $ (\pi')^{-1}(b)$.
\end{enumerate}
\br \label{rem how to define m} The above proof  shows that $m_b$ divides $\mu_b$.  To explain this fact, we fist recall the proof of \cite[Proposition 1.4]{ACM10}. 
Consider the following commutative diagram
$$
\xymatrix{
 \pi_1(f^{-1}(U)) \ar[r] \ar[d] &               \pi_1(U) \ar[d] \\
  \pi_1(X)\ar[r] &   \pi_1(C).
 }
 $$
All the vertical maps are surjective since they are induced by inclusions.
All the horizontal maps are surjective since the generic fiber is connected.
For any point $b\in B$, say  $f^{*}(b)\coloneqq  \mu^1_b E_1 +\cdots +\mu^k_b E_k $, where $\{E_1,\cdots,E_k\}$ are reduced irreducible components of $f^{-1}(b)$.  Then the point $b$ has the signed multiplicity $\mu_b\coloneqq \gcd(\mu^1_b, \cdots, \mu_b^k)$.
Since $f$ maps $\gamma_{E_i}$ to $(\gamma_b)^{\mu^i_b}$, where $\gamma_{E^i}$ and $\gamma_b$ are the meridians for $E_i$ and $b$, respectively,  it is easy to see that the first horizontal map in the above diagram factors through $\pi_1(X)\twoheadrightarrow\pi_1^{\orb}(C) $.

By the same proof, we get a surjective map $\pi_1(T_b)\twoheadrightarrow \Z/\mu_b \Z$. In particular, we get the following commutative diagram
$$
\xymatrix{
\pi_1(T_b)  \ar[r] \ar[d] &             \pi_1(X)   \ar[d] & \\
  \Z/\mu_b \Z\ar[r] &   \pi_1^\orb(C) \ar[r] & H,
 }
 $$
 where $\Z/\mu_b \Z \to   \pi_1^\orb(C) $ is given by the natural  injective map. Hence $m_b$ can be defined as the order of the image group for the composition of the bottom horizontal maps. In particular, $m_b$ divides $\mu_b$.
\er
\begin{proof}[Proof of Theorem \ref{main thm}(b)]
To finish the proof, we need to compute $\chi(S)$. Note that
\be
\chi(S)-\sum_{b\in B}\frac{|T|}{m_b}=|T|(\chi(C)-\vert B \vert),
\ee
hence
\be
\chi(S)=|T|(\chi(C)+\sum_{b\in B}(\frac{1}{m_b}-1)).
\ee
Using the results in Section 2 and Lemma \ref{lem finite cover}, we get
\be \label{formula}
\alpha_1(X^{\nu},\K)=-\chi(C)+\sum_{b\in B}(1-\frac{1}{m_b})+\sum_{b\in B, p\mid\frac{\mu_b}{m_b}}\frac{1}{m_b}.
\ee
and
\be
\M_1(X^{\nu},\K)=\sum_{b\in B}\frac{1}{m_b}\log\frac{\mu_b}{m_b}.
\ee
\end{proof}

\br \label{rem average} Let $\K$ be an algebraically closed filed with characteristic $p\geq 0$.
If $p$ does not divide $\vert T \vert$, then $R\pi_* \K_{X^T}$ is a direct sum of  $\vert T\vert$-many rank one local systems and $\Hom(H,\K^*)$ has exactly $\vert T \vert$ many connected components, where each of them is a copy of $\Hom(\Z^n ,\K^*) $. By \cite[Theorem 4.8]{LL21}, the formula (\ref{formula}) should be understood as the average of  $\dim H^1(X, \nu^{-1} L_\lambda)$ for $\lambda$ being general in every connected components.
\er
\bex{\cite[Example 6.12]{DS14}} \label{ex key example}
Fix an integer $\mu\geq 2$. Let $\sA_\mu$ be the deleted monomial arrangement,
where its defining equations in $\mathbb{CP}^2$ are $yz(x^\mu-y^\mu)(x^\mu-z^\mu)(y^\mu-z^\mu)  $. Ordering the
hyperplanes as the factors of the defining polynomial. Its projective complement $X$ admits an orbifold map $X\to \C^*$ given by
$$\dfrac{z^\mu(x^\mu -y^\mu)}{y^\mu(x^\mu-z^\mu)}, $$
which is of type $(0,2,\mu)$. 

Let $m$ be a positive integer, which divides $\mu$. Consider $\nu $ as the composition of the following maps
$$\pi_1(X)\twoheadrightarrow \Z*\Z/\mu \Z \twoheadrightarrow \Z  \oplus \Z/m \Z ,$$
where the first map is induced by the orbifold map and the second map is the quotient map.
Let $\K$ be an algebraically closed field with characteristic $p\geq 0$.
Then we have \begin{center}
$\alpha_1(X^\nu,\K)=\begin{cases}
1, & \mathrm{if}\ p\mid \frac{\mu}{m},\\
 1-\frac{1}{m}, & \mathrm{ otherwise}
\end{cases}  $
\end{center}
and $\M_1(X^\nu)=  \frac{1}{m}\log\frac{\mu}{m}$.

On the other hand, if $p$ does not divide $m$, then \cite[Theorem 4.8]{LL21} shows that for general $\lambda \in \Hom(\Z\oplus \Z/m \Z, \K^*)$,
\begin{center}
$\dim H_1(X,\nu^{-1}L_\lambda)=\begin{cases}
0, & \mathrm{if}\ p \nmid \mu \text{ and } \lambda \text{ is general in }  \Hom(\Z, \K^*),\\
1, &  \mathrm{ otherwise},
\end{cases}  $
\end{center}
 where its average is same as the formula given above.
\eex

\end{document}